\documentclass[12pt, reqno]{amsart}
\usepackage{amssymb}
\usepackage{amsfonts}
\usepackage{amssymb,amsmath,amscd}
\usepackage[colorlinks=true]{hyperref}
\usepackage{mathrsfs}
\newtheorem{Theorem}{\indent Theorem}[section]

\newtheorem{Lemma}[Theorem]{\indent Lemma}

\theoremstyle{remark}
\newtheorem{Remark}{Remark}

\begin{document}
\centerline{\bf On  fractional sums  of the divisor functions
}

\bigskip
\centerline{Wei Zhang}

\centerline{
School of Mathematics and Statistics,
Henan University,
Kaifeng  475004, Henan, China}
\centerline{zhangweimath@126.com}
\bigskip

\textbf{Abstract}\
In this paper, we consider  the fractional sum of the divisor functions.
We can improve previous results  considered by Bordell\`{e}s \cite{Bo}   and Liu-Wu-Yang \cite{LWY}.
Precisely, we can show that
 \begin{align*}
 S_{\tau_{k}}(x)=\sum_{n\leq x}\tau_{k}\left(\left[\frac{x}{n}\right]\right)=
\sum_{n=1}^{\infty}
\frac{\tau_{k}(n)}{n(n+1)}x
+O\left(x^{9/19+\varepsilon}\right),
 \end{align*}
 where $\varepsilon$ is an arbitrary small positive constant
and $\tau_{k}(n)$ is the number of representations of $n$ as product of $k$ natural numbers.



\bigskip
\bigskip
\numberwithin{equation}{section}

\section{Introduction}
Recently, the sum
\[
S_{f}(x)=\sum_{n\leq x}f\left(\left[\frac{x}{n}\right]\right)
\]
has attracted  many experts special attention (for example, see \cite{BDHPS,Bo,LWY,Wu,Zhai}),
where $f$ is a complex-valued arithmetic function and $[\cdot]$ denotes the floor function (i.e. the greatest integer function).
 One can call $S_{f}(x)$ the fractional sum of $f$ (see \cite{St}).

Specially, for some fixed $\eta\in(0,1)$ and
\[
f(n)\ll n^{\eta},
\]
independently, Wu \cite{Wu} and Zhai \cite{Zhai} showed that
\[
S_{f}(x)=C_{f}x+O\left(x^{(1+\eta)/2}\right),
\]
where $f$ is a complex-valued arithmetic function and
\[
C_{f}=\sum_{n=1}^{\infty}\frac{f(n)}{n(n+1)}.
\]
This formula improves the recent result obtained by
Bordell\`{e}s, Dai, Heyman, Pan and Shparlinski \cite{BDHPS}.
In fact, for general function $f$  being a positive real-valued arithmetic function, then one can obtain some much better results by involving the theory of Fourier series \cite{Va} and exponential sums \cite{GK}. For example, one can refer to  \cite{Bo,LWY,MS,St}.

Let $\varepsilon$ be  an arbitrarily small positive constant. And this  positive constant  may be different for different situations in this article.
By using the symmetry of the divisor function, in \cite{MS},  it is proved that
\begin{align*}
 S_{\tau }(x)=\sum_{n\leq x}\tau \left(\left[\frac{x}{n}\right]\right)=
\sum_{n=1}^{\infty}
\frac{\tau (n)}{n(n+1)}x
+O\left(x^{11/23+\varepsilon}\right),
 \end{align*}
where $\tau (n)$ is  the number of representations of $n$ as product of two natural numbers and
\[
11/23\approx0.4782.
\]
Recently, this result was improved by many experts.
In \cite{Bo}, by using the Dirichlet hyperbolic method and more effort, it is proved that
\begin{align*}
 S_{\tau}(x)=\sum_{n\leq x}\tau\left(\left[\frac{x}{n}\right]\right)=
\sum_{n=1}^{\infty}
\frac{\tau(n)}{n(n+1)}x
+O\left(x^{19/40+\varepsilon}\right),
 \end{align*}
 where
$$19/40=0.475.$$
By using a new estimate on 3-dimensional exponential sums, in \cite{LWY}, Liu-Wu-Yang showed that
\begin{align*}
 S_{\tau}(x)=\sum_{n\leq x}\tau\left(\left[\frac{x}{n}\right]\right)=
\sum_{n=1}^{\infty}
\frac{\tau(n)}{n(n+1)}x
+O\left(x^{9/19+\varepsilon}\right),
 \end{align*}
 where
$$9/19\approx0.47368.$$
This sum has relations  to the generalized divisor function $\tau_{k}(n)$ ($k\geq2$), where  $\tau_{k}(n)$ is the number of representations of $n$ as product of $k$ natural numbers.
Hence, in \cite{Bo} and \cite {LWY}, some results for these arithmetic functions are also given.
Precisely, in \cite{LWY}, it is proved that
%\begin{align}\label{tauk}
% S_{2^{\omega}}(x)=\sum_{n\leq x}2^{\omega\left(\left[x/n\right]\right)}=
%\sum_{n=1}^{\infty}
%\frac{2^{\omega\left(n\right)}}{n(n+1)}x
%+O\left(x^{9/19+\varepsilon}\right)
% \end{align}
% and
 \begin{align}\label{2w}
 S_{\tau_{k}}(x)=\sum_{n\leq x}\tau_{k}\left(\left[\frac{x}{n}\right]\right)=
\sum_{n=1}^{\infty}
\frac{\tau_{k}(n)}{n(n+1)}x
+O\left(x^{\theta(k)+\varepsilon}\right),
 \end{align}
 where
 \[
 \theta(k)=\frac{5k-1}{10k-1}.
 \]
 These results improve previous results of \cite{Bo} and can be seen as a generalization of $S_{\tau}(x).$
More recently, by using some deep  results of Jutila \cite{Ju}, Stucky \cite{St} showed that
\begin{align}\label{tau2}
 S_{\tau}(x)=\sum_{n\leq x}\tau\left(\left[\frac{x}{n}\right]\right)=
\sum_{n=1}^{\infty}
\frac{\tau(n)}{n(n+1)}x
+O\left(x^{5/11+\varepsilon}\right),
 \end{align}
 where
$$5/11\approx0.454545.$$
This improves previous results of $S_{\tau}(x).$
It is reasonable to believe that one can obtain much better results for  $S_{\tau_{k}}(x)$  by using the method of Stucky \cite{St}. However, the method of Stucky \cite{St} is too special to be generalized to the   situation of
$S_{\tau_{k}}(x)$.

The purpose of this paper is to study $S_{\tau_{k}}(x)$.
We can improve   (\ref{2w}) by giving the  following   result, which gives   improvement for
 \[
 \theta(k)=\frac{5k-1}{10k-1}.
 \]
\begin{Theorem}\label{th2}
Let $\tau_{k}(n)$ be the number of representations of $n$ as product of $k$ natural numbers.
Then we have
 \begin{align*}
 S_{\tau_{k}}(x)=\sum_{n\leq x}\tau_{k}\left(\left[\frac{x}{n}\right]\right)=
\sum_{n=1}^{\infty}
\frac{\tau_{k}(n)}{n(n+1)}x
+O\left(x^{9/19+\varepsilon}\right),
 \end{align*}
 where $\varepsilon$ is an arbitrary small positive constant.
\end{Theorem}
\begin{Remark}
In fact, a little better result can also be given. One can refer to \cite{LWY1} for details. For any exponent pair   $(\kappa,\lambda),$ one can obtain that
 \begin{align*}
 S_{\tau_{k}}(x)=\sum_{n\leq x}\tau_{k}\left(\left[\frac{x}{n}\right]\right)=
\sum_{n=1}^{\infty}
\frac{\tau_{k}(n)}{n(n+1)}x
+O\left(x^{\frac{2\kappa+\lambda+3}{4\kappa+\lambda+7}+\varepsilon}\right).
\end{align*}
If we choose $(\kappa,\lambda)=(1/2,1/2),$ we can obtain Theorem \ref{th2}. If we choose $(\kappa,\lambda)=(1653/3494+\varepsilon,
1760/3494+\varepsilon)
=BA^{5}(13/84+\varepsilon,55/84+\varepsilon),$
we can obtain a little better result. This basic observation can also be seen in \cite{LWY1}.
If we assume that $(\varepsilon,1/2+\varepsilon)$ is also an exponent pair, then we can obtain
 \begin{align*}
 S_{\tau_{k}}(x)=\sum_{n\leq x}\tau_{k}\left(\left[\frac{x}{n}\right]\right)=
\sum_{n=1}^{\infty}
\frac{\tau_{k}(n)}{n(n+1)}x
+O\left(x^{\frac{7}{15}+\varepsilon}\right).
\end{align*}
\end{Remark}

\section{Proof of Theorem \ref{th2}}
We will start  the proof for Theorem \ref{th2}  with some necessary lemmas. The following lemma is from \cite{LWY}.
\begin{Lemma}\label{LWY}
Let $\alpha>0,$ $\beta>0,$ $\gamma>0$ and $\delta\in\mathbb{R}$ be some constants. For $X>0,$ $H\geq1,$ $M\geq 1,$ and $N\geq 1,$ define
\[
S_{\delta}=S_{\delta}(H,M,N):=
\sum_{h\sim H}\sum_{m\sim M}\sum_{n\sim N}
a_{h,n}b_{m}e\left(X\frac{M^{\beta}N^{\gamma}}{H^{\alpha}}
\frac{h^{\alpha}}{m^{\beta}n^{\gamma}+\delta}\right),
\]
where $e(t)=e^{2\pi i t},$ the $a_{h,n}$ and $b_{m}$ are complex  numbers such that $a_{h,n}\leq 1,$ $b_{m}\leq 1$ and $m\sim M$ means that $M<m\leq 2M.$ For any $\varepsilon>0,$ we have
\[
S_{\delta}\ll
\left(\left(X^{\kappa}
H^{2+\kappa}M^{1+\kappa+\lambda}
N^{2+\kappa}\right)^{1/(2+2\kappa)}
+HM^{1/2}N+H^{1/2}MN^{1/2}+X^{-1/2}HMN\right)
X^{\varepsilon}
\]
uniformly for $M\geq 1,$ $N\geq 1,$ $H\leq N^{\gamma-1}M^{\beta}$ and $|\delta|\leq 1/\varepsilon,$ where
$(\kappa,\lambda)$ is an exponent pair and the implied constant may depend on $\alpha,\beta,\gamma,$ and $\varepsilon.$
\end{Lemma}
Next lemma can be seen in Theorem A.6 in \cite{GK} or Theorem 18 in \cite{Va}.
\begin{Lemma}\label{z3}
For $0<|t|<1,$ let
$$W(t) = \pi t(1-|t|)\cot\pi t + |t|.$$ For $x\in\mathbb{R},$ $H\geq1,$ we define
$$\psi^{*}(x)=\sum_{1\leq |h|\leq H}(2\pi ih)^{-1}W\left(\frac{h}{H+1}\right)e(hx)$$
and
\[
\delta(x)=\frac{1}{2H+2}\sum_{|h|\leq H}\left(1-\frac{|h|}{H+1}\right)e(hx).
\]
Then $\delta(x)$ is non-negative, and we have
$$|\psi^{*}(x)-\psi(x)|\leq
\delta(x).$$
\end{Lemma}
 
We also need
the following  well-known lemma (for example, one can refer to page 441 of \cite {BB} or page 34 of  \cite{GK}).
\begin{Lemma}\label{z2}
Let $g^{(m')}(x)\asymp YX^{1-m'}$ for $1< X\leq x\leq 2X$ and $m'=1,2,\cdots.$ Then one has
\[
\sum_{X<n\leq 2X}e(g(n))\ll Y^{\kappa}X^{\lambda}+Y^{-1},
\]
where $(\kappa,\lambda)$ is any exponent pair.
\end{Lemma}

Now we begin the proof of Theorem \ref{th2}. Let
\[
N=x^{9/19}.
\]
We can write
\[
S_{\tau_{k}}(x):=S_{\tau_{k},1} +S_{\tau_{k},2} ,
\]
where
\begin{align}\label{c11}
S_{\tau_{k},1}=\sum_{n\leq N}\tau_{k}\left(\left[\frac{x}{n}\right]\right)
\end{align}
and
\begin{align}\label{c12}
S_{\tau_{k},2}=\sum_{N<n\leq x}\tau_{k}\left(\left[\frac{x}{n}\right]\right).
\end{align}
Obviously, by   $\tau_{k}(n)\ll n^{\varepsilon},$ we have
\begin{align*}
S_{\tau_{k},1}=\sum_{n\leq N}\tau_{k}\left(\left[\frac{x}{n}\right]\right)
&=\sum_{n\leq N}(x/n)^{\varepsilon}
\\
&\ll N^{1+\varepsilon}\\&\ll x^{9/19+\varepsilon}.
\end{align*}

As to $S_{\tau_{k},2},$  by    $\tau_{k}(n)\ll n^{\varepsilon},$
  we have
$$
\sum_{n\leq x} \tau_{k}(n)  \ll x^{1+\varepsilon}.
$$
Hence we can get
\begin{align}\label{c1}
\begin{split}
S_{\tau_{k},2}&=\sum_{N<n\leq x}\tau_{k}\left(\left[\frac{x}{n}\right]\right)
\\&
=\sum_{d\leq x/N}\tau_{k}(d)\sum_{x/(d+1)<n\leq x/d}1
\\&
=\sum_{d\leq x/N}\tau_{k}(d)\left(\frac{x}{d}-\frac{x}{d+1}
-\psi(\frac{x}{d})
+\psi(\frac{x}{d+1})\right)
\\&=x\sum_{d=1}^{\infty}\frac{\tau_{k}(d)}{d(d+1)}
+O\left(
N^{1+\varepsilon}
\right)+(\log x)\sum_{D<d\leq 2D}\tau_{k}(d)\psi\left(\frac{x}{d+\delta}\right),
\end{split}\end{align}
where $N\leq D\leq x/N$ and $\delta\in\{0,1\}.$
Then we need to estimate
\[\sum_{D<d\leq 2D}\tau_{k}(d)\psi\left(\frac{x}{d+\delta}\right).
\]
By Lemma \ref{z3}, we have
\begin{align}\label{c2}
\begin{split}\sum_{D<d\leq 2D}&\tau_{k}(d)\psi\left(\frac{x}{d+\delta}\right)\\
&\ll
\sum_{1\leq h\leq H}\frac1h\sum_{D<d\leq 2D}\tau_{k}(d)e\left(\frac{hx}{d+\delta}\right)
\\&+ \sum_{1\leq h\leq H}\frac1H\sum_{D<d\leq 2D}\tau_{k}(d)e\left(\frac{hx}{d+\delta}\right)
+D/H.
\end{split}\end{align}
Then we will focus on the estimate of
\begin{align*}
S_{\delta}:=\sum_{1\leq h\leq H}\frac1h\sum_{D<d\leq 2D}\tau_{k}(d)e\left(\frac{hx}{d+\delta}\right).
\end{align*}
By using the relation
$$\sum_{n_{1}n_{2}\cdots n_{k}=n}1=\tau_{k}(n),$$
and the dichotomy  method, we have
\[
S_{\delta}\ll D^{\varepsilon}\sum_{1\leq h\leq H}\frac1h\sum_{d_{i}\sim D_{i},i=1,2,\cdots, k}e\left(\frac{hx}{d_{1}d_{2}\cdots d_{k}+\delta}\right),
\]
where
\begin{align}\label{k1}
d_{i}\leq d_{i+1}, \ D_{i}\leq D_{i+1}, \ \textup{for} \  1\leq i \leq k-1\end{align}
and
\begin{align}
\label{k2}
\prod_{i=1}^{k}D_{i}\sim D.\end{align}

Now we divide three cases to deal with this.

{\bf Case I}

Suppose that $D_{k}\geq D^{2/3}.$

By Lemma \ref{z2} and choosing $(\kappa,\lambda)=(1/2,1/2)$, we have
\begin{align*}
S_{\delta}&\ll D^{\varepsilon}\sum_{1\leq h\leq H}\frac1h\sum_{d_{i}\sim D_{i},i=1,2,\cdots, k-1}
\left(\left(\frac{hx}{d_{1}d_{2}\cdots d_{k-1}D_{k}^{2}}\right)^{1/2}
D_{k}^{1/2}\right.
\\
&\left.+\frac{(d_{1}d_{2}\cdots d_{k-1}D_{k})^{2}}{hx}\right)
\\&\ll D^{2/9+\varepsilon}x^{1/3}
+D^{2+\varepsilon}/x,
\end{align*}
where  we have chosen   $H=D^{7/9}x^{-1/3}.$
We assume that $D>x^{1/2}.$ Hence we can verify that $H>1.$
Then for  $D_{k}\geq D^{2/3},$ by choosing $N=x^{9/19},$ we have
\[
S_{\delta} \ll D^{2/9}x^{1/3}\ll x^{77/171+\varepsilon}\ll x^{9/19}.
\]

{\bf Case II}

Suppose that  $D^{1/3}\leq D_{k}\leq D^{2/3}.$

By choosing $N=x^{9/19}$ and  $(\kappa,\lambda)=(1/2,1/2)$ in Lemma \ref{LWY}, and restricted the range to $D^{1/3}\leq D_{k}\leq D^{1/2}$ according to symmetry, we have
\begin{align*}
S_{\delta}&\ll
D^{\varepsilon}\sum_{1\leq h\leq H}\frac1h\sum_{d_{i}\sim D_{i},i=1,2,\cdots, k-1}\sum_{d_{k}\sim D_{k}}e\left(\frac{hx}{d_{1}d_{2}\cdots d_{k}+\delta}\right) \\
&\ll
(DH)^{\varepsilon}
\left(x^{1/6}D^{1/2}D^{1/12}+D^{3/4}
+x^{-1/2}D^{3/2}\right)
\\
&\ll x^{9/19 +\varepsilon}.
\end{align*}

 {\bf Case III}

 Suppose that $D_{k}\leq D^{1/3}.$
 Then by (\ref{k1}) and (\ref{k2}), we have $D_{i}\leq D^{1/3},$ $i=1,2,\ldots,k.$ We   also suppose that $t$ is the least integer such that $D_{1}D_{2}\ldots D_{t}>D^{1/3}.$ Then we have
 \[
 D^{1/3}\leq (D_{1}D_{2}\ldots D_{t-1})D_{t}\leq D^{2/3}.
 \]
 Let $l_{1}=d_{1}d_{2}\ldots d_{t}$ and let $l_{2}=d_{t+1}d_{t+2}\ldots d_{k}.$
 Then we have
\[
S_{\delta}\ll
D^{\varepsilon}\sum_{1\leq h\leq H}\frac1h\sum_{l_{1}\sim L_{1} }d_{t}(l_{1})\sum_{l_{2}\sim L_{2} }d_{k-t}(l_{2})e\left(\frac{hx}{l_{1}l_{2} +\delta}\right),
\]
where $D^{1/3}\leq L_{1}\leq D^{2/3}$
 and $D^{1/3}\leq L_{2}\leq D^{2/3}.$
Then similar as the second case, we have
\begin{align*}
S_{\delta}&\ll
(DH)^{\varepsilon}
\left(x^{1/6}D^{1/2}D^{1/12}+D^{3/4}
+x^{-1/2}D^{3/2}\right)
\\ & \ll x^{9/19 +\varepsilon}.
\end{align*}

Then from the above three cases, we have
\[S_{\delta}:=\sum_{1\leq h\leq H}\frac1h\sum_{D<d\leq 2D}\tau_{k}(d)e\left(\frac{hx}{d+\delta}\right)\ll x^{9/19+\varepsilon}.
\]
Then by (\ref{c1})-(\ref{c2}), we have
\[S_{\tau_{k}}(x)= x\sum_{d=1}^{\infty}\frac{\tau_{k}(d)}{d(d+1)}
+ O\left(x^{9/19+\varepsilon}\right).
\]
Recall   (\ref{c11}) and (\ref{c12}), then we can finally give Theorem \ref{th2}.

\bigskip
\bigskip
{\bf Acknowledgements} The author would like to thank the the referee who gives some  detailed corrections and suggestions. The author also would like to thank Professor J. Wu for sending some interesting related papers. This work was supported by NSFC(11871307).

\address{Wei Zhang\\ School of Mathematics and Statistics\\
               Henan University\\
               Kaifeng  475004, Henan\\
               China}
\email{zhangweimath@126.com}

\end{document}